\documentclass[11pt]{article}
\usepackage{latexsym}
\usepackage{epsfig}
\usepackage{amssymb, amsmath, multicol}
\usepackage[dvips]{color}
\usepackage{graphicx}
\newtheorem{lemma}{Lemma}
\newtheorem{proposition}{Proposition}
\newtheorem{theorem}{Theorem}
\renewcommand{\theequation}{\arabic{equation}}

\def\bbeta{\mbox{\boldmath{$\beta$}}}
\def\toD{\overset{\text{D}}{\longrightarrow}}

\newcommand{\PE}{\mbox{PE}}

\newcommand{\Var}{\mbox{Var}}
\newcommand{\Cov}{\mbox{Cov}}
\newcommand{\cor}{\mbox{Cor}}

\newcommand{\etal}{{\em et al.\,}\/}
\newcommand{\bY}{\mbox{\bf Y}}

\title{\bf\Large Dynamic Integration of Time- and State-domain Methods for Volatility
Estimation}

\author{
By JIANQING FAN\\
Benheim  Center for Finance and Department of ORFE\\
Princeton University,
Princeton, NJ 08544\\
jqfan@princeton.edu
\and
 YINGYING FAN\\
Department of ORFE,
Princeton University,
Princeton, NJ 08544\\
yingying@princeton.edu
\and
AND JIANCHENG JIANG
\\
LMAM and School of
Mathematical Sciences,
Peking University, Beijing 100871\\
Jiang@math.pku.edu.cn
}

\date{}
\begin{document}
\maketitle

\begin{center} Summary \end{center}

Time- and state-domain methods are two common approaches for
nonparametric prediction. The former predominantly uses the data
in the recent history while the latter mainly relies on historical
information. The question of combining these two pieces of
valuable information is an interesting challenge in statistics.
We surmount this problem
via dynamically integrating information from both the time and the state domains.
The estimators from both domains are optimally combined based on a
data driven weighting strategy, which provides a more efficient
estimator of volatility.
Asymptotic normality is seperately
established for the time damain, the state domain, and the
integrated estimators. By comparing the efficiency of the
estimators, it is demonstrated that the proposed integrated
estimator uniformly dominates the two other estimators. The
proposed dynamic integration approach is also applicable to other
estimation problems in time series. Extensive simulations are
conducted to demonstrate that the newly proposed procedure
outperforms some popular ones such as the RiskMetrics and the
historical simulation approaches, among others. Empirical studies
endorse convincingly our integration method.\\

\noindent{\it Some key words:} Bayes; Dynamical integration;
State-domain; Time-domain; Volatility.

\newpage
\oddsidemargin 0.0 in

\begin{center}
\section{Introduction}
\end{center}

In forecasting a future event or making an investment decision,
two pieces of useful information are frequently consulted. Based
on the recent history, one uses a form of local average, such as
the moving average in the time-domain, to forecast a future event.
This approach uses the continuity of a function and ignores
completely the information in the remote history, which is related
to current through stationarity. On the other hand, one can
forecast a future event based on state-domain modeling such as the
ARMA, TAR, ARCH models or nonparametric models (see Tong, 1990;
Fan \& Yao, 2003 for details). For example, to forecast the
volatility of the yields of a bond with the current rate 6.47\%,
one computes the standard deviation based on the historical
information with yields around 6.47\%. This approach relies on the
stationarity and depends completely on historical data. But, it
ignores the importance of the recent data. The question of how to
combine the estimators from both the time-domain and the
state-domain poses an interesting challenge to statisticians.

\begin{figure}[htbp]
\centerline{\psfig{figure=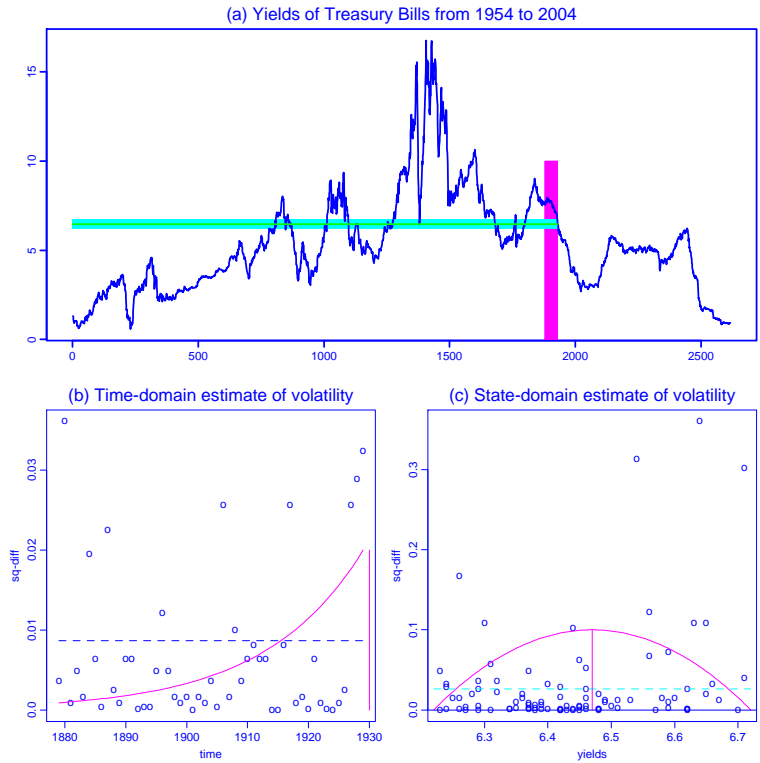,width=4in}}
\begin{singlespace}
\caption[Fig1] {\small Illustration of time and state-domain
estimation. (a)  The yields of 3-month treasury bills from 1954 to
2004.  The vertical bar indicates localization in time and the
horizontal bar represents localization in the state. (b)
Illustration of time-domain smoothing:  squared differences are
plotted against its time index and the exponential weights are
used to compute the local average.  (c)  Illustration of the
state-domain smoothing:  squared differences are plotted against
the level of interest rates, restricted to the interval $6.47\%
\pm .25\%$ indicated by the horizontal bar in Figure 1(a). The
Epanechnikov kernel is used for computing the local average. }
\label{fig1}
\end{singlespace}

\end{figure}

To elucidate our idea, consider the weekly data on the yields of
3-month treasury bills presented in Figure 1.  Suppose that the
current time is January 04, 1991 and interest rate is 6.47\% on
that day, corresponding to the time index $t = 1930$.  One may
estimate the volatility based on the weighted squared differences
in the past 52 weeks (1 year), say.  This corresponds to the
time-domain smoothing, using a small vertical stretch of data in
Figure 1(a). Figure 1(b) computes the squared differences of the
past year's data and depicts its associated exponential weights.
The estimated volatility (conditional variance) is indicated by
the dashed horizontal bar.   Let the resulting estimator be
$\hat{\sigma}_{t, \mbox{\scriptsize time}}^2$. On the other hand,
in financial activities, we do consult historical information in
making better decisions.  The current interest rate is 6.47\%. One
may examine the volatility of the yields when the interest rates
are around 6.47\%, say, $6.47\% \pm .25\%$.  This corresponds to
using the part of data indicated by the horizontal bar.  Figure
1(c) plots the squared differences $X_t - X_{t-1}$ against
$X_{t-1}$ with $X_{t-1}$ restricted to the interval $6.47\% \pm
.25\%$. Applying the local kernel weight to the squared
differences results in a state-domain estimator $\hat{\sigma}_{t,
\mbox{\scriptsize state}}^2$, indicated by the horizontal bar in
Figure 1(c). Clearly, as shown in Figure 1(a),  except in the
3-week period right before January 4, 1991 (which can be excluded
in the state domain fitting), the last period with interest rate
around $6.47\% \pm .25\%$ is the period from May 15, 1988 and July
22, 1988.  Hence, the time and state-domain estimators use two
nearly independent components of the time series, as they are
136-week apart in time. See the horizontal and vertical bars of
Figure 1(a). These two kinds of estimators have been used in the
literature for forecasting volatility.   The former is prominently
featured in the RiskMetrics of J.P. Morgan, and the latter has
been used in nonparametric regression (see Tong, 1995; Fan \& Yao,
2003 and references therein).  The question arises how to
integrate them.

An integrated estimator is to introduce a dynamic weighting scheme
$0 \leq w_t \leq 1$ to combine the two nearly independent
estimators. Define the resulting integrated estimators as
$$
\hat{\sigma}_t^2 = w_t \hat{\sigma}_{t, \mbox{\scriptsize time}}^2
+ (1-w_t) \hat{\sigma}_{t, \mbox{\scriptsize state}}^2.
$$
The question is how to choose the dynamic weight $w_t$ to optimize
the performance.  A reasonable approach is to minimize the
variance of the combined estimator, leading to the dynamic optimal
weights
\begin{equation}
w_t = \frac{\Var(\hat{\sigma}_{t, \mbox{\scriptsize state}}^2)}{
\Var(\hat{\sigma}_{t, \mbox{\scriptsize time}}^2) +
\Var(\hat{\sigma}_{t, \mbox{\scriptsize state}}^2)}, \label{a1}
\end{equation}
since the two piece of estimators are nearly independent. The
unknown variances in (\ref{a1}) can easily be estimated in Section
3. Another approach is the Bayesian approach, which regards the
historical information as the prior. We will explore this idea in
Section 4. The proposed method is also applicable to other
estimation problems in time series such as forecasting the mean
function and the volatility matrix of multivariate time series.

To appreciate the intuition behind our approach, let us consider
the diffusion process
\begin{equation}
     dr_t = \mu(r_t) dt + \sigma(r_t) dW_t, \label{a2}
\end{equation}
where $W_t$ is a Wiener process. This diffusion process is
frequently used to model asset price and the yields of bonds,
which are fundamental to fixed income securities, financial
markets, consumer spending, corporate earnings, asset pricing and
inflation. The family of models include famous ones such as the
Vasicek (1977) model, the CIR model (Cox, \etal  1985) and the CKLS
model (Chan, \etal  1992). Suppose that at time $t$ we have a historic data
$\{r_{t_i}\}_{i =0}^{N}$ from the process (\ref{a2}) with a
sampling interval $\Delta$. Our aim is to estimate the volatility $\sigma_t^2\equiv \sigma^2(r_t).$
Let $Y_i =\Delta^{-1/2} (r_{t_{i+1}} - r_{t_i})$. Then for the model (\ref{a2}),
the Euler approximation scheme is
\begin{equation}
Y_i \approx \mu(r_{t_i}) \Delta^{1/2} + \sigma(r_{t_i}) \varepsilon_i,
\label{a3}
\end{equation}
where $\varepsilon_i \sim_{i.i.d.} N(0, 1)$ for $i=0,\cdots,N-1$.  Fan \& Zhang (2003)
studied the impact of the order of difference on statistical
estimation.  They found that while higher order can possibly
reduce approximation errors, it increases variances of data
substantially.  They recommended the Euler scheme (\ref{a3}) for
most practical situations.  The time-domain smoothing relies on
the smoothness of $\sigma(r_{t_i})$ as a function of time $t_i$.  This
leads to the exponential smoothing estimator in Section 2.1.  On
the other hand, the state-domain smoothing relies on structural
invariability implied by the stationarity:  the conditional
variance of $Y_i$ given $r_{t_i}$ remains the same even for the data
in the history.  In other words, historical data also furnish the
information about $\sigma(\cdot)$ at the current time. Combining
these two nearly independent estimators leads to a better
estimator.

In this paper, we focus on the estimation of volatility of a
portfolio to illustrate how to deal with the problem of dynamic
integration.  Asymptotic normality of the proposed estimator is
established and extensive simulations are conducted, which
theoretically and empirically demonstrate the dominated
performance of the integrated estimation.\\


{\centering \section{Estimation of Volatility}}

The volatility estimation is an important issue of modern
financial analysis since it pervades almost every facet of this
field. It is a measure of risk of a portfolio and is related to
the Value-at-Risk (VaR), asset pricing, portfolio allocation,
capital requirement and risk adjusted returns, among others. There
is a large literature on estimating the volatility based on
time-domain and state-domain smoothing.   For an overview, see the
recent book by Fan \& Yao (2003). \\

{\centering \subsection{Time-domain estimator}}

A popular version of time-domain estimator of the volatility is
the moving average estimator:
\begin{equation}
\hat{\sigma}_{MA, t}^2 = n^{-1} \sum_{i=t-n}^{t-1} Y_i^2,
\label{b1}
\end{equation}
where $n$ is the size of the moving window. This estimator ignores
the drift component, which contributes to the variance in the
order of $O(\Delta)$ instead of $O(\Delta^{1/2})$  (see Stanton,
1997 and Fan \& Zhang, 2003), and utilizes local $n$ data points.
An extension of the moving average estimator is the exponential
smoothing estimation of the volatility given by
\begin{equation}
  \hat{\sigma}_{ES, t}^2 = (1 - \lambda) Y_{t-1}^2 + \lambda
  \hat{\sigma}_{ES, t-1}^2
= (1 - \lambda) \{ Y_{t-1}^2 + \lambda Y_{t-2}^2 +
 \lambda^2 Y_{t-3}^2 +
 \cdots \}, \label{b2}
\end{equation}
where $\lambda$  is a smoothing parameter  that controls the size
of the local neighborhood. The RiskMetrics of J.P. Morgan (1996),
which is used for measuring the risks, called Value at Risk (VaR),
of financial assets, recommends $\lambda = 0.94$ and $\lambda =
0.97$ respectively for calculating VaR of the daily and monthly
returns.

The exponential smoothing estimator in (\ref{b2}) is a weighted
sum of the squared returns prior to time $t$. Since the weight
decays exponentially, it essentially uses recent data. A slightly
modified version that explicitly uses only $n$ data points before
time $t$ is
\begin{equation}
  \hat{\sigma}_{ES,t}^2 =
  \frac{1-\lambda}{1-\lambda^n}\sum_{i=1}^nY_{t-i}^2\lambda^{i-1}.
  \label{b3}
\end{equation}
When $\lambda = 1$, it becomes the moving average estimator
(\ref{a1}).  With slight abuse of notation, we will also denote
the estimator for $\sigma^2(r_t)$ as $\hat{\sigma}_{ES,t}^2$.

All of the time domain smoothing is based on the assumption that
the returns $Y_{t-1},\ Y_{t-2},$ $\cdots,\ Y_{t-n}$ have
approximately the same volatility.  In other words, $\sigma(r_t)$
in (\ref{a1}) is continuous in time $t$.  The following
proposition gives the condition under which this holds.

\begin{proposition}\label{P1}  Under Conditions (A1) and (A2) in the
Appendix, we have
\begin{equation*}
|\sigma^2(r_s)-\sigma^2(r_u)|\leq K |s-u|^{(p-1)/(2p)},
\end{equation*}
for any $s,u \in [t-\eta, t]$, where the coefficient $K$ satisfies
$E[K^{2(p+\delta)}]<\infty$ and $\eta$ is a positive constant.
\end{proposition}

With the above H\"{o}lder continuity, we can establish the
asymptotic normality of the time-domain estimator.

\begin{theorem}\label{T1}
Suppose that $\sigma^2_t>0$. Under conditions (A1) and (A2), if
$n\rightarrow +\infty$  and $n\Delta\rightarrow 0$, then
$$
  \hat{\sigma}_{ES,t}^2 - \sigma^2_t{\longrightarrow} 0,\, \
 \mbox{\rm a.e.}
 $$
Moreover, if  the limit $c=\lim_{n\to \infty}n(1-\lambda)$ exists
and $n\Delta^{(p-1)/(2p-1)} \rightarrow 0$,
$$
   \sqrt{n}[\hat{\sigma}_{ES,t}^2 - \sigma^2_t]/ s_{1,t}
   \stackrel{\mathcal D}{\longrightarrow}
   \mathcal{N}\left(0,1\right),
$$
where $s_{1,t}^2 =c\,\sigma^4_t\frac{e^c+1}{e^c-1}.$
\end{theorem}

Theorem~\ref{T1} has very interesting implications.  Even though
the data in the local time-window is highly correlated (indeed,
the correlation tending to one), we can compute the variance as if
the data were independent.  Indeed, if the data in (\ref{b3}) were
independent and locally homogeneous, we have
\begin{eqnarray*}
\Var( \hat{\sigma}_{ES,t}^2) & \approx &
\frac{(1-\lambda)^2}{(1-\lambda^n)^2} 2 \sigma^4_t
\sum_{i=1}^n \lambda^{2(i-1)} \\
& = & \frac{2 \sigma^4_t (1 - \lambda) (1 + \lambda^n)}{
      (1+\lambda) (1 - \lambda^n)}
\approx  \frac{ 1}{n} s_{1,t}^2.
\end{eqnarray*}
This is indeed the asymptotic variance given in Theorem~\ref{T1}.

\begin{center}
\subsection{Estimation in state-domain}\label{non}
\end{center}

To obtain the nonparametric estimation of the functions $f(x)=
\Delta^{1/2} \mu(x)$ and $\sigma^2(x)$ in (\ref{a3}), we use the
local linear smoother studied in Ruppert \etal  (1997) and  Fan
\& Yao (1998). The local linear technique is chosen for its
several nice properties, such as the asymptotic minimax efficiency
and the design adaptation. Further, it automatically corrects edge
effects and facilitates the bandwidth selection (Fan \& Yao,
2003).

To facilitate the theoretical argument in Section 3, we exclude
the $n$ data points used in the time-domain fitting.   Thus, the
historical data at time $t$ are $\{(r_{t_i}, Y_i), i = 0, \cdots,
N-n-1\}$.  Let $\hat f(x)=\hat{\alpha}_1$ be the local linear
estimator that solves the following weighted least-squares
problem:
\[
(\hat{\alpha}_1,\hat{\alpha}_2)=\arg\min_{\alpha_1,\alpha_2}\sum_{i=0}^{N-n-1}
[Y_{i}-\alpha_1-\alpha_2(r_{t_i}-x)]^2K_{h_1}(r_{t_i}-x), \] where
$K(\cdot)$ is a kernel function and $h_1>0$ is a bandwidth. Denote
the squared residuals by $\hat R_i= \{Y_i-\hat f(r_{t_i})\}^2$. Then
the local linear estimator of $\sigma^2(x)$ is $\hat{\sigma}_{S}^2
(x) = \hat{\beta}_0$ given by
\begin{equation}
(\hat{\beta}_0,\hat{\beta}_1)=\arg\min_{\alpha,\
\beta}\sum_{i=0}^{N-n-1}
\{\hat{R}_{i}-\beta_0-\beta_1(r_{t_i}-x)\}^2W_h(r_{t_i}-x)
 \label{b4}
\end{equation}
with kernel function $W$ and  bandwidth $h$.  Fan \& Yao (1998)
gives strategies of bandwidth selection.  It was shown in Stanton
(1997) and Fan \& Zhang (2003) that $Y_i^2$ instead of
$\hat{R}_i$ in (\ref{b4}) can also be used for the estimation of
$\sigma^2 (x)$.


The asymptotic bias and variance of $\hat{\sigma}_{S}^2 (x)$ are
given by Fan \& Zhang (2003, theorem 4).  Set $\nu_j=\int
u^jW^2(u)du$ for $j=0,1,2$. Let $p(\cdot)$ the invariant density
function of the Markov process $\{r_s\}$ from (\ref{a1}).  Then,
we have

\begin{theorem}\label{T2}
Let $x$ be in the interior of the support of $p(\cdot)$. Suppose
that the second derivatives $\mu(\cdot)$ and $\sigma^2(\cdot)$
exist in a neighborhood of $x$. Under conditions (A3)-(A7), we
have
\begin{align*}
\sqrt{(N-n)h}[\hat{\sigma}^2_{S} (x) - \sigma^2(x)] / s_{2}(x)
\stackrel{\mathcal D}{\longrightarrow} {\mathcal N} \left(0, 1
\right) ,
\end{align*}
where $s_2^2(x) =2\nu_0\sigma^4(x)/p(x).$
\end{theorem}

\begin{center}
\section{Dynamic Integration of time and state domain estimators}\label{sec3}
\end{center}

In this section, we first show how the optimal dynamic weights in (\ref{a1}) can be
estimated and then prove that the time-domain and state-domain
estimator are indeed asymptotically independent.\\

{\centering\subsection{Estimation of dynamic weights}}

For the exponential smoothing estimator in (\ref{b3}), we can
apply the asymptotic formula given in Theorem~\ref{T1} to get an
estimate of its asymptotic variance.  However, since the estimator
is a weighted average of $Y_{t-i}^2$, we can obtain its variance
directly by assuming $Y_{t-j} \sim N(0, \sigma_t^2)$ for small
$j$. Indeed, with the above local homogeneous model, we have
\begin{eqnarray}
   \Var(\hat{\sigma}_{ES, t}^2)
& \approx & \frac{(1-\lambda)^2}{(1-\lambda^n)^2} 2 \sigma_t^4
     \sum_{i=1}^{n} \sum_{j=1}^{n} \lambda^{i+j-2} \rho(|i-j|) \nonumber \\
& = & \frac{2 (1-\lambda)^2  \sigma_t^4 }{(1-\lambda^n)^2} \{ 1 +
2 \sum_{k=1}^{n-1} \rho(k) \lambda^k (1 - \lambda^{2(n-k)})
   /(1-\lambda^2) \},  \label{c1}
\end{eqnarray}
where $\rho(j) = \cor(Y_t^2, Y_{t-j}^2)$ is  the autocorrelation
of the series $\{Y_{t-j}^2\}$.  The autocorrelation can be
estimated from the data in history.  Note that due to the locality
of the exponential smoothing, only $\rho(j)$'s with the first 30
lags, say, contribute to the variance calculation.

We now turn to estimate the variance of $\hat{\sigma}_{S, t}^2=
\hat{\sigma}_S^2(r_{t})$. Details can be found in Fan \& Yao
(1998) and \S6.2 of Fan \& Yao (2003). Let
$$
V_j(x) = \sum_{i=1}^{t-1} (r_{t_i} -x)^j
W\Bigl(\frac{r_{t_i}-x}{h_1}\Bigr)
$$
and
$$
    \quad \xi_i(x) = W\Bigl(\frac{r_{t_i}-x}{h_1}\Bigr) \{V_2(x) -
        (r_{t_i}-x) V_1(x) \}/\{ V_0(x) V_2(x) - V_1(x)^2\}.
$$
Then the local linear estimator can be expressed as
$$
   \hat{\sigma}_S^2(x) = \sum_{i=1}^{t-1} \xi_i(x) \hat{R}_i
$$
and its variance can be approximated as
\begin{equation}
     \Var(\hat{\sigma}_S^2(x)) \approx \Var \{ (Y_1 - f(x))^2 | r_{t_1}=x \}
     \sum_{i=1}^{t-1} \xi_i^2(x). \label{c2}
\end{equation}
See also Figure 1 and the discussions at the end of \S2.1. Again,
for simplicity, we assume that $\Var(\hat R_i|r_{t_i}=x) \approx 2
\sigma^4(x)$, which holds if $\varepsilon_t \sim N(0,1)$.

Combining (\ref{a1}), (\ref{c1}) and (\ref{c2}), we propose to
combine the time-domain and the state-domain estimator with the
dynamic weight
\begin{equation}
\hat{w}_t = \frac{ \hat{\sigma}_{S,t}^4 \sum_{i=1}^{t-1}
\xi_i^2(r_{t})} {\hat{\sigma}_{S,t}^4
 \sum_{i=1}^{t-1} \xi_i^2(r_t)
+ c_t \hat{\sigma}_{ES, t}^4
}, \label{c3}
\end{equation}
where $c_t = \frac{(1-\lambda)^2 }{(1-\lambda^n)^2} \{ 1 + 2
\sum_{k=1}^{n-1} \rho(k) \lambda^k (1 - \lambda^{2(n-k)})
/(1-\lambda^2) \}$ [see (\ref{c1})]. This is obtained by
substituting (\ref{c1}) and (\ref{c2})  into (\ref{a1}). For
practical implementation, we truncate the series
$\{\rho(i)\}_{i=1}^{t-1}$ in the summation as
$\{\rho(i)\}_{i=1}^{30}$.   This results in the dynamically
integrated estimator
\begin{equation}
\hat{\sigma}_{I, t}^2 = \hat{w}_t \hat{\sigma}_{ES,t}^2 + (1 -
\hat{w}_t) \hat{\sigma}_{S,t}^2,  \label{c4}
\end{equation}
where $\hat{\sigma}_{S,t}^2 = \hat{\sigma}_S^2(r_t)$. The function
$\hat{\sigma}_S^2(\cdot)$ depends on the time $t$ and we need to
update this function as time evolves. Fortunately, we need only to
know the function at the point $r_t$. This reduces significantly
the computational cost. The computational cost can be reduced
further, if we update the estimated function
$\hat{\sigma}_{S,t}^2$ at a prescribed time schedule (e.g. once
every two months for weekly data).


Finally, we would like to note that in the choice of weight, only
the variance of the estimated volatility is considered, rather
than the mean square error.  This is mainly to facilitate the
dynamically weighted procedure.  Since the smoothing parameters in
$\hat{\sigma}_{ES, t}^2$ and $\hat{\sigma}^2_S(x)$ have been tuned
to optimize their performance separately, their biases and
variances trade-off have been considered.  Hence, controlling the
variance of the integrated estimator $\hat{\sigma}_{I,t}^2$ has
also controlled, to some extent, the bias of the estimator. Our
method focuses only on the estimation of volatility, but the
method can be adapted to other estimation problems, such as the
value at risk studied in Duffie \& Pan (1997) and the drift
estimation for diffusion considered in Spokoiny (2000) and
volatility matrix for multivariate time series. Further study
along this topic is beyond the scope of the current investigation.

{\centering\subsection{Sampling properties}}

The fundamental component to the choice of dynamic weights is the
asymptotic independent between the time and state-domain
estimator.  By ignoring the drift term (see Stanton, 1997; Fan \&
Zhang 2003), both the estimators $\hat{\sigma}_{ES, t}^2$ and
$\hat{\sigma}_{S,t}^2$ are linear in  $\{Y_i^2\}$.  The following
theorem shows that the time-domain and state-domain estimators are
indeed asymptotically independent.  To facilitate the notation, we
present the result at the current time $t_N$.

\begin{theorem}\label{T3}
Let $s_{2,t_N}=s_2(r_{t_N}).$ Under the conditions of Theorems
\ref{T1} and \ref{T2}, if the condition (A2) holds at point $t_N$,
we have
\begin{itemize}
\item [(a)] asymptotic independence:
$$
[\sqrt{n}(\hat{\sigma}^2_{ES, t_N} -
\sigma^2_{t_N})/s_{1,t_N},\sqrt{(N-n)h}(\hat{\sigma}^2_{S,t_N} -
\sigma^2_{t_N})/s_{2,t_N}]^T \stackrel{\mathcal
D}{\longrightarrow} {\mathcal N}(0,I_2).
$$

\item [(b)] asymptotic normality of $\hat{\sigma}_{I,t_N}^2$:
 if the limit $d=\lim_{N\to\infty}n/[(N-n)h]$ exists, then
$$
   \sqrt{(N-n)h / \omega} [\hat{\sigma}^2_{I, t_N}-\sigma^2_{t_N})]
   \stackrel{\mathcal D}{\longrightarrow} {\mathcal N}(0,1),
$$
where $\omega=w_{t_N}^2 s_{1,t_N}^2/d+(1-w_{t_N})^2 s_{2,t_N}^2$.
\end{itemize}
\end{theorem}






From Theorem \ref{T3}, based on the optimal weight the asymptotic
relative efficiencies of $\hat{\sigma}_{I,t_N}^2$ with respect to
$\hat{\sigma}_{S,t_N}^2$ and $\hat{\sigma}_{ES,t_N}^2$ are
respectively
$$
\mbox{eff}(\hat{\sigma}_{I,t_N}^2,\hat{\sigma}_{S,t_N}^2)=1 + d
s_{2,t_N}^2/s_{1,t_N}^2, \qquad \mbox{and} \qquad
\mbox{eff}(\hat{\sigma}_{I,t_N}^2,\hat{\sigma}_{ES,t_N}^2)= 1 +
s_{1,t_N}^2/(ds_{2,t_N}^2),
$$
which are greater than one. This demonstrates that the integrated
estimator $\hat{\sigma}_{I,t_N}^2$ is more efficient than the time
domain and the state domain estimators.

\begin{center}
\section{Bayesian integration of volatility estiamtes}\label{sec4}
\end{center}

Another possible approach is to consider the historical information
as the prior and to incorporate them in the estimation of
volatility by the Bayesian framework.  We now explore such an
approach.

{\centering\subsection{Bayesian estimation of volatility}
\label{bay} }

The Bayesian approach is to regard the recent data $Y_{t-n},
\cdots, Y_{t-1}$ as an independent sample from $N(0,\sigma^2)$
[see (\ref{a3})] and to regard the historical information being
summarized in a prior. To incorporate historical information, we
assume that the variance $\sigma^2$ follows an Inverse Gamma
distribution with parameters $a$ and $b$, which has the density
function
$$
f(\sigma^2)=b^a\Gamma^{-1}(a){\{\sigma^2\}}^{-(a+1)}\mbox{exp}
(-b/\sigma^2).
$$
Denote by $\sigma^2 \sim IG(a, b)$.  It is a well-known fact that
\begin{equation}
\mbox{E}(\sigma^2)=\frac{b}{(a-1)},~~~
\mbox{Var}(\sigma^2)=\frac{b^2}{(a-1)^2(a-2)},~~~
\mbox{mode}(\sigma^2)=\frac{b}{(a+1)}. \label{d1}
\end{equation}
The hyperparameters $a$ and $b$ will be estimated from historical
data such as the state-domain estimators.

It can easily be shown that the posterior density of $\sigma^2$
given $\bY=(Y_{t-n}, \cdots, Y_{t-1})$ is
IG$(a^*,b^*)$, where
$$
    a^*=a+\frac {n}{2}, \quad b^*=\frac {1}{2}\sum_{i=1}^n Y_{t-i}^2+b.
$$
From (\ref{d1}), the Bayesian mean of $\sigma^2$ is
$$
\hat\sigma^2 = \frac{b^*}{(a^*-1)}=
   \sum_{i=1}^n (Y_{t-i}^2+2b)/(2(a-1)+n).
$$
This Bayesian estimator can easily be written as
\begin{equation}
\hat{\sigma}^2_B = \frac{n}{n+2(a-1)} \hat{\sigma}_{MA, t}^2
        + \frac{2(a-1)}{n + 2(a-1)} \hat{\sigma}_P^2,
        \label{d2}
\end{equation}
where $\hat{\sigma}_{MA, t}^2$ is the moving average estimator
given by (\ref{b1}) and $\hat{\sigma}_P^2 = b/(a-1)$ is the prior
mean, which will be determined from the historical data. This
combines the estimate based on the data and prior knowledge.

The Bayesian estimator (\ref{d3}) utilizes the local average of
$n$ data points. To incorporate the exponential smoothing
estimator (\ref{b2}), we regard it as the local average of
\begin{equation}
   n^*= \sum_{i=1}^n \lambda^{i-1} = \frac{1 - \lambda^n}{1 - \lambda}
    \label{d3}
\end{equation}
data points. This leads to the following integrated estimator
\begin{eqnarray}
\hat{\sigma}_{B, t}^2 & = &
    \frac{n^*}{n^*+2(a - 1)}\hat{\sigma}_{ES,t}^2 +
        \frac{2(a - 1)}{2(a -1)+n^*} \hat{\sigma}^2_{P} \nonumber \\
     &=& \frac{1-\lambda^n}{1 - \lambda^n + 2(a - 1)(1-\lambda)}\hat{\sigma}_{ES,t}^2
        + \frac{2(a - 1)(1 - \lambda)}
        {1 - \lambda^n + 2(a - 1)(1 - \lambda)}\hat{\sigma}_{P}^2.
        \label{d4}
\end{eqnarray}
In particular, when $\lambda = 1$, the estimator (\ref{d4})
reduces to (\ref{d2}).\\

{\centering\subsection{Estimation of Prior Parameters} }

A reasonable source for obtaining the prior information in
(\ref{d4}) is based on the historical data up to time $t$. Hence,
the hyper-parameters $a$ and $b$ should depend on $t$ and can be
used to match with the historical information. Using the
approximation model (\ref{a3}), we have
\[
E[(Y_t- \hat f(r_t))^2\mid r_t] \approx  \sigma^2(r_t) \qquad
\Var[(Y_t- \hat f(r_t))^2\mid r_t] \approx 2 \sigma^4(r_t).
\]
These can be estimated from the historical data up to time $t$,
namely, the state-domain estimator $\hat{\sigma}^2_S(r_t)$. Since we
have assumed that prior distribution for $\sigma_t^2$ is IG($a_t,
b_t)$, then by the method of moments, we would set
\begin{eqnarray*}
  E( \sigma_t^2 )=  \frac{b_t}{a_t-1} = \hat \sigma_S^2(r_t),\\
  \Var(\sigma_t^2)= \frac{b_t^2}{(a_t-1)^2 (a_t-2)}
   = 2 \hat{\sigma}_S^4 (r_t).
\end{eqnarray*}
Solving the above equation, we obtain that
$$
  \hat{a}_t= 2.5\, \,\ \ \mbox{\rm and}\, \,\ \ \hat{b}_t=1.5\hat \sigma_S^2(r_t).
$$
Substituting this into (\ref{d4}), we obtain the following
estimator
\begin{equation}
\hat{\sigma}_{B, t}^2 = \frac{1-\lambda^n}{1 - \lambda^n + 3
    (1-\lambda)}\hat{\sigma}_{ES,t}^2
        + \frac{3(1 - \lambda)}
        {1 - \lambda^n + 3(1 - \lambda)}\hat{\sigma}_{S,t}^2.
        \label{d5}
\end{equation}
Unfortunately, the weights in (\ref{d5}) are static, which does not
depend on the time $t$.  Hence, the Bayesian method does not
produce a satisfactory answer to this problem.\\

{\centering\section{Numerical Analysis }\label{sec5} }

To facilitate the presentation, we use the simple abbreviation in
Table \ref{tab1} to denote five volatility estimation methods.
Details of the first three methods can be found in Fan \& Gu
(2003).  In particular, the first method is to estimate the
volatility using the standard deviation of the yields in the past
year and the RiskMetrics method is based on the exponential
smoothing with $\lambda = 0.94$.  The semiparametric method of Fan
\& Gu (2003) is an extension of a local model used in the
exponential smoothing, with the smoothing parameter determined by
minimizing the prediction error. It includes the exponential
smoothing with $\lambda$ selected by data as a specific example.

\begin{table}[htbp]
\begin{center}
\caption{Abbreviations of five volatility estimators  \label{tab1}} \doublerulesep
0.5pt \small
\begin{tabular}{@{}ll@{}} \hline \hline
Hist: the historical method \\
RiskM: the RiskMetrics method of J.P. Morgan \\
Semi: the semiparametric estimator (SEV) in Fan \& Gu (2003)\\
NonBay: the nonparametric Bayesian method in (\ref{d5}) with $\lambda=0.94$ \\
Integ: the integration method of time and state domains in (\ref{c4}) \\
 \hline
\hline
\end{tabular}
\end{center}
\end{table}

The following four measures are employed to assess the performance
of different procedures for estimating the volatility. Other
related measures can also be used. See Dav\'e \& Stahl (1997).

{\bf Measure 1}. Exceedence ratio against confidence level.

This measure counts  the number of the events for which the loss
of an asset exceeds the loss predicted by the normal model at a
given confidence $\alpha$.  With estimated volatility, under the
normal model, the one-period VaR is estimated by $\Phi^{-1}
(\alpha) \hat{\sigma}_{t}$, where $\Phi^{-1}(\alpha)$ is the
$\alpha$ quantile of the standard normal distribution. For each
estimated VaR,  the Exceedence Ratio (ER) is computed as
\begin{equation}
\mbox{ER}(\hat{\sigma}^2_t) = m^{-1} \sum_{i=T+1}^{T+m}
      I(Y_{i} <  \Phi^{-1} (\alpha) \hat{\sigma}_{i}), \label{er1}
\end{equation}
for an out-sample of size $m$. This gives an indication on how
effective the volatility estimator can be used for predicting the
one-period VaR. Note that the Monte Carlo error for this measure
has an approximate size $\{\alpha (1-\alpha)/m \}^{1/2}$, even
when the true $\sigma_t$ is used.  For example, with $\alpha =
5\%$ and $m = 1000$, the Monte Carlo error is around $0.68\%$.
Thus, unless the post-sample size $m$ is large enough, this
measure has difficulty in differentiating the performance of
various estimators due to the presence of large error margins.
Note that the ER depends strongly on the assumption of normality.
If the underlying return process is non-normal, the Student's
$t(5)$ say, the ER will grossly be overestimated even with the
true volatility.
In our simulation study, we will employ the true $\alpha$-quantile
of the error distribution instead of $\Phi^{-1}(\alpha)$ in
(\ref{er1}) to compute the ER.  For real data analysis, we use the
$\alpha$-quantile of the last $250$ residuals for the in-sample
data.

{\bf Measure 2}. Mean Absolute Deviation Error.

To motivate this measure, let us first consider the mean square
errors:
$$
  \PE(\hat{\sigma}^2_t) =  m^{-1} \sum_{i=T+1}^{T+ m}  (Y_{i}^2 -
  \hat{\sigma}_{i}^2)^2.
$$
The expected value can be decomposed as
\begin{equation}
E (\PE) =  m^{-1} \sum_{i=T+1}^{T+ m}  E ( \sigma_{i}^2 -
  \hat{\sigma}_i^2)^2 + m^{-1} \sum_{i=T+1}^{T+ m} E ( Y_{i}^2 -
  \sigma_{i}^2)^2.\label{b15}
\end{equation}
Note that the first term reflects the effectiveness of the
estimated volatility while the second term is the size of the
stochastic error, independent of estimators. As in all statistical
prediction problems, the second term is usually of an order of
magnitude larger than the first term.  Thus, a small improvement
on PE could mean substantial improvement over the estimated
volatility.  However, due to the well-known fact that financial
time series contain outliers, the mean-square error is not a
robust measure.  Therefore, we used the mean-absolute deviation
error (MADE):
$$
\mbox{MADE}(\hat{\sigma}^2_t) = m^{-1} \sum_{i=T+1}^{T+ m} \mid
Y_i^2 - \hat \sigma_i^2
                   \mid.
$$

{\bf Measure 3.} Square-root Absolute Deviation Error.

An alternative variation to MADE is the square-Root Absolute
Deviation Error (RADE), which is defined as
$$
\mbox{RADE}(\hat{\sigma}^2_t) = m^{-1}  \sum_{i=T+1}^{T+m} \Bigl |
             \mid Y_i\mid- \sqrt{\frac{2}{\pi}}\hat\sigma_i \Bigr |.
$$
The constant factor comes from the fact that $E|\varepsilon_t| =
\sqrt{\frac{2}{\pi}}$ for $\varepsilon_t \sim N(0,1)$. If the
underlying error distribution deviates from normality, this
measure is not robust.

{\bf Measure 4.} Ideal Mean Absolute Deviation Error.

To assess the estimation of the volatility in simulations, one can
also employ the ideal mean absolute deviation error (IMADE):
$$
\mbox{IMADE}=m^{-1}\sum_{i=T+1}^{T+m}
|\hat{\sigma}^2_i-\sigma^2_i|.
$$
This measure calibrates the accuracy of the forecasted volatility
in terms of the absolute difference between the true and the
forecasted one. However, for real data analysis, this measure is
not applicable.

{\centering \subsection{Simulations}}

To assess the performance of the five estimation methods in Table
\ref{tab1}, we compute the average and the standard deviation of
each of the four measures over $600$ simulations. Generally
speaking, the smaller the average (or the standard deviation), the
better the estimation approach. We also compute the ``score'' of
an estimator, which is the percentage of times among 600
simulations that the estimator outperforms the average of the 5
methods in terms of an effectiveness measure. To be more specific,
for example, consider RiskMetrics using MADE as an effectiveness
measure. Let $m_i$ be the MADE of the RiskMetrics estimator at the
$i$-th simulation, and $\bar{m}_i$ the average of the MADEs for
the five estimators at the $i$-th simulation. Then the ``score''
of the RiskMetrics approach in terms of the MADE is defined as
$$\frac{1}{600}\sum_{i=1}^{600} I(m_i< \bar{m}_i).$$
Obviously, the estimators with higher scores are preferred. In
addition, we define a ``relative loss'' of an estimator
$\hat{\sigma}^2_t$ relative to $\hat{\sigma}^2_{I,t}$
 in terms of MADEs as
$$\mbox{\rm RLOSS}(\hat{\sigma}^2_t\,\hat{\sigma}^2_{I,t})
=\frac{\overline{\mbox{MADE}}(\hat{\sigma}^2_t)-
   \overline{\mbox{MADE}}(\hat{\sigma}^2_{I,t})}
      {\overline{\mbox{MADE}}(\hat{\sigma}^2_{I,t})},$$
where $\overline{\mbox{MADE}}(\hat{\sigma}^2_t)$ is the average of
MADE($\hat{\sigma}^2_t$) among simulations.

{\bf Example 1.} To simulate the interest rate data, we consider
the Cox-Ingersoll-Ross (CIR) model:
\begin{eqnarray*}
dr_t=\kappa(\theta-r_t)dt+\sigma r_t^{1/2}dW_t, \ \ t\geq t_0,
\end{eqnarray*}
where the spot  rate, $r_t$, moves around a central location or
long-run equilibrium level $\theta=0.08571$ at speed
$\kappa=0.21459$. The $\sigma$ is set to be 0.07830. These values
of parameters are cited from Chapman \& Pearson (2000), which
satisfy the condition $2\kappa\theta\ge \sigma^2$ so that the
process  $r_t$ is stationary and positive.
The model has been studied by Chapman \& Pearson (2000)
and Fan \& Zhang (2003).

There are two methods to generate samples from this model. The first one
is the discrete-time order $1.0$ strong approximation scheme in
Kloeden, \etal (1996);  the
second one is using the exact transition density detailed in Cox et al. (1985) and
 Fan \& Zhang (2003). Here we use the first method to generate
$600$ series of data each with length $1200$ of the weekly
data from this model. For each simulation, we set the first $900$ observations
as the ``in-sample'' data
and the last $300$ observations as the ``out-sample'' data.

\begin{table}[htbp]
\begin{center}
\caption{Comparisons of several volatility estimation methods \label{tab2}}
\doublerulesep 0.5pt \small
\begin{tabular}{@{}| c| c| c|c|c|c|c|@{}} \hline \hline
 Measure & Empirical Formula & Hist & RiskM & Semi & NonBay &Integ\\
\hline\hline
       & Score (\%)             & 17.17   & 20.83   & 32.00  & 44.33  & 99.83\\
 IMADE    & Ave $(\times 10^{-5})$ & 0.2383 & 0.2088  & 0.1922 & 0.1833 & 0.0879\\
       & Std $(\times 10^{-5} )$ & 0.1087 & 0.0746  & 0.0718 & 0.0675 & 0.0554 \\
       & Relative Loss (\%)      & 171.20 & 137.61  & 118.79 & 108.60  & 0 \\
\hline
       & Score (\%)             &  39.83  & 54.33   & 60.00  & 57.17  & 72.17\\
MADE   & Ave $(\times 10^{-4})$ &  0.1012 & 0.0930  & 0.0932 & 0.0924 & 0.0903\\
       & Std$(\times 10^{-5})$  &  0.3231 & 0.3152  & 0.3010 & 0.3119 & 0.2995\\
       & Relative Loss (\%)     &  12.03  & 2.95    & 3.16   & 2.31   & 0\\
\hline
       & Score (\%)             &  40.83  & 53.33   & 54.83  & 57.50  & 74.50\\
RADE   & Ave                    & 0.0015  & 0.0015  & 0.0015 & 0.0015 & 0.0014  \\
       & Std ($\times 10^{-3}$)  & 0.2530  & 0.2552  & 0.2461 & 0.2536 & 0.2476\\
       & Relative Loss (\%)     & 6.88    & 1.66    & 2.13   & 1.27   & 0\\
\hline
ER     & Ave                    & 0.0556  & 0.0547  & 0.0536 & 0.0535  & 0.0508 \\
       & Std                    & 0.0257  & 0.0106  & 0.0122 & 0.0107  & 0.0122 \\
\hline \hline
\end{tabular}
\end{center}
\end{table}

\begin{figure}[htbp]
\centerline{\epsfig{figure=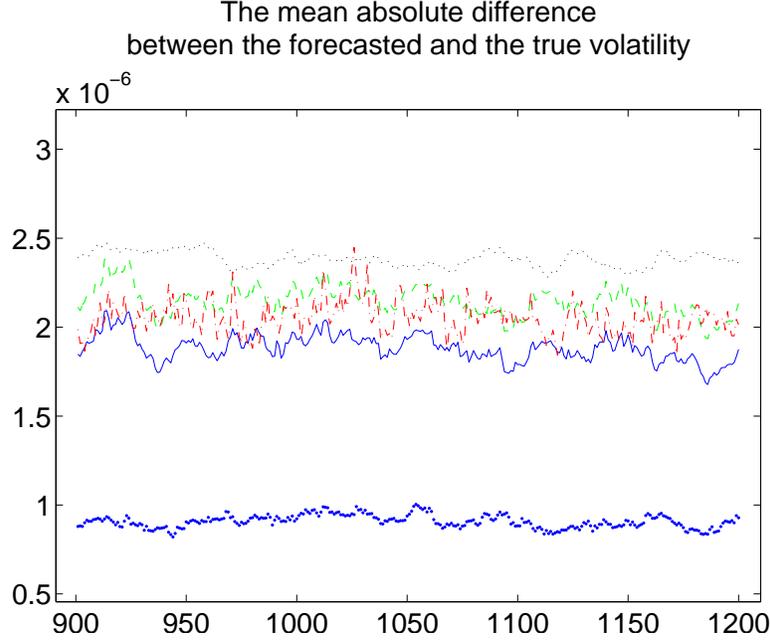,width=4in}}
\begin{singlespace}
\caption[Fig2] {The mean absolute difference between the
forecasted and the true volatility.  
Solid - integrated estimator (11); 
small circle - nonparametric Bayesian integrated estimator (16); 
star - historical method; 
dashed - RiskMetrics; 
dash dotted - Semiparametric estimator in Fan \& Gu
(2003).} \label{fig2}
\end{singlespace}
\end{figure}

The results are summarized in Table \ref{tab2}, which shows that
the performance of the integrated estimator uniformly dominates
the other estimators because of its highest score, lowest IMADE,
MADE, and RADE.   The improvement in IMADE is over
$100$ percent. This shows that our integrated volatility method
better captures the volatility dynamics. The Bayesian method of
combining the estimates from the time and state domains outperforms all
other methods.   The historical simulation
method performed poorly due to mis-specification of the function
of the volatility parameter. The results here show the advantage
of aggregating the information of time domain and state domain.
Note that all estimators have reasonable ER values at level
$0.05$, especially the ER value of the integrated estimator is
closest to $0.05$.
To appreciate how much improvement for our integrated method over the other methods,
we display the mean absolute difference between the forecasted and the true volatility in Figure \ref{fig2}.
It is seen that the integrated method is much better than the others in terms of the difference.

{\bf Example 2}. There is a large literature on the estimation of
volatility.  In addition to the famous parametric models such as
ARCH and GARCH, stochastic volatility models have also received a
lot of attention.  For an overview, see, for example,
Barndoff-Neilsen \& Shephard (2001, 2002), Bollerslev \& Zhou
(2002) and references therein.  We consider the following
stochastic volatility model:
\begin{eqnarray*}
&&dr_t=\sigma_t dB_t, \ r_0=0 \\
&&dV_t=\kappa(\theta-V_t)dt+\alpha V_tdW_t, \ V_0=\eta,\
V_t=\sigma_t^2,
\end{eqnarray*}
where
$W_t$ and $B_t$ are two independent standard Brownian motions.

There are two methods to generate samples from this model. One is
the direct method, using the result of  Genon-Catalot et al.
(1999). Let $a=1+2\kappa/\alpha^2$ and $b=2\theta\kappa/\alpha^2$.
The conditions (A1)-(A4) in the above paper are satisfied with the
parameter values in the model being constants as $\kappa=3$,
$\theta=0.009$ and $\alpha^2=4$ and the initial random variable
$\eta$ follows the Inverse Gamma distribution.
The value of $\theta$ is set as the real variance of the daily
return for Standard \& Poor 500 data from January 4, 1988 to
December 29, 2000. The value $\alpha^2$ is to make the parameter
$a$ of the stable distribution $IG(a,b)$ equal $2.5$, the prior
parameter in (\ref{c3}).  If $\Delta \rightarrow 0$ and $n\Delta
\rightarrow \infty$,  then
$$Y_i\rightarrow \sqrt{\frac{b}a}T,\ \mbox{where}\  T\sim t(2a).$$

Another method is the discretization of the model.
Conditionally on $\textsl{g}=\sigma(V_t,t\geq 0)$, the random variables
$Y_i$ are independent and follow $N(0,\bar{V}_i)$ with
$$\bar{V}_i=\frac 1 \Delta \int_{(i-1) \Delta}^{i \Delta}V_s ds.$$
To simulate the diffusion
process $V_t$, one can use the following order 1.0 scheme with
sampling interval $\Delta^*=\Delta/30$,
\begin{eqnarray*}
V_{i+\Delta^*}=V_{i} + \kappa(\theta-V_{i})\Delta^* + \alpha V_{i}
(\Delta^*)^{1/2}
    \varepsilon_i +
    \frac{1}{2} \alpha^2 V_{i} \Delta^* (\varepsilon_i^2 -1),
\end{eqnarray*}
where $\{\varepsilon_i\}$ are independent random series from the
standard normal distribution.

We simulate $600$ series of $1000$ monthly data using the
second method with step size $\Delta = 1/12$. For each simulated
series, set the first three quarters observations as the in-sample
data and the remaining observations as the out-sample data. The
performance of each volatility estimation is described in Table
\ref{tab3}. The conclusion similar to Example 1 can be drawn from
this example.

\begin{table}[htbp]
\begin{center}
\caption{Comparisons of several volatility estimation methods\label{tab3}}
\doublerulesep 0.5pt \small
\begin{tabular}{@{}|c| c|c|c|c|c|  c|@{}} \hline \hline
 Measure & Empirical Formula & Hist & RiskM & Semi & NonBay &Integ\\
\hline\hline
       & Score (\%)             & 27.67   & 49.33  & 52.83  & 58.83  & 77.17\\
IMADE    & Ave                    & 0.0056  & 0.0051 & 0.0051 & 0.0050 & 0.0047\\
       & Std                    & 0.0023  & 0.0019 & 0.0021 & 0.0018 & 0.0016\\
       & Relative Loss (\%)     & 17.74   & 7.63   & 6.56   & 5.18   & 0\\
\hline
       & Score (\%)             & 35.33   & 52.17  & 57.67  & 58.00  & 82.67\\
MADE   & Ave                    & 0.0099  & 0.0089 & 0.0087 & 0.0088 & 0.0082 \\
       & Std                    & 0.0032  & 0.0022 & 0.0022 & 0.0021 & 0.0017 \\
       & Relative Loss (\%)     & 20.48   & 7.53   & 5.38   & 6.17   & 0 \\
\hline
       & Score (\%)             & 33.00   & 49.17  & 53.33  & 58.83  & 81.33\\
RADE   & Ave                    & 0.0477  & 0.0455 & 0.0452 & 0.0451 & 0.0438 \\
       & Std                    & 0.0059  & 0.0051 & 0.0051 & 0.0049 & 0.0042\\
       & Relative Loss (\%)     & 8.77    & 3.70   & 3.11   & 2.91   & 0\\
\hline
ER     & Ave                    & 0.0457  & 0.0547 & 0.0546 & 0.0516 & 0.0533\\
       & Std                    & 0.0156  & 0.0126 & 0.0143  & 0.0127 & 0.0146\\
\hline \hline
\end{tabular}
\end{center}
\end{table}


{\bf Example 3.} 
We now
consider the geometric Brownian (GBM):
$$dr_t=\mu r_t+\sigma r_t dW_t,$$
where $W_t$ is a standard one-dimensional Brownian motion.
This is a non-stationary process to which we check if our method
continues to apply.
Note that the celebrated Black-Scholes option price formula is derived
on the Osborne's assumption that the stock price follows the GBM
model. By the It$\hat{o}$ formula, we have
$$\log r_t-\log r_0=(\mu-\sigma^2/2)t+\sigma^2 W_t.$$
We set $\mu=0.03$ and $\sigma=0.26$ in our simulations. With the
Brownian motion simulated from independent Gaussian increments,
one can generate the samples for the GBM.
Here we use the
latter with $\Delta=1/52$ in $600$ simulations. For each
simulation, we generate $1000$ observations and use the first two
thirds of observations as in-sample data and the remaining
observations as out-sample data.

\begin{table}[htbp]
\begin{center}
\caption{Comparisons of several volatility estimation methods \label{tab4}}
\doublerulesep 0.5pt \small
\begin{tabular}{@{}|c| c|c|c|c|c|c|@{}} \hline \hline
 Measure & Empirical Formula & Hist & RiskM & Semi & NonBay &Integ\\
\hline\hline
       & Score (\%)             & 2.17   & 89.98  & 7.01 & 99.17   & 99.17\\
 IMADE    & Ave $(\times 10^{-5})$ & 0.1615 & 0.0811 & 0.1154 & 0.0746 &  0.0746\\
       & Std $(\times 10^{-4} )$ &  0.1030 & 0.0473 & 0.0632 & 0.0440 &  0.0440\\
       & Relative Loss (\%)      & 116.42  & 8.64  & 54.63 &  0   & 0\\
\hline
       & Score (\%)             &   40.17  & 58.67 & 54.00 &  60.00    & 66.17\\
MADE   & Ave $(\times 10^{-5})$ &  0.2424  & 0.2984 & 0.2896 & 0.2958 &   0.2859\\
       & Std$(\times 10^{-4})$  & 0.1037  & 0.1739 & 0.1633 & 0.1723  &  0.1663\\
       & Relative Loss (\%)     & -15.24  & 4.35 &  1.30 & 3.46 & 0\\
\hline
       & Score (\%)             &  36.83   & 60.17  & 47.50  & 62.33  & 69.50\\
RADE   & Ave  ($\times 10^{-3}$)  &  0.5236  & 0.4997 & 0.5114 & 0.4975 & 0.4903\\
       & Std ($\times 10^{-3}$) & 0.5898  &  0.6608 &  0.6567  &  0.6573  &   0.6435 \\
      & Relative Loss (\%)     & 6.80   & 1.92  & 4.30  & 1.47   & 0\\
\hline
ER     & Ave                    &  0.0693  & 0.0532 &  0.0517 &  0.0506  &  0.0444\\
       & Std                    &  0.0467  & 0.0095 &  0.0219 &  0.0110  &   0.0160\\
\hline \hline
\end{tabular}
\end{center}
\end{table}

\begin{table}[htbp]
\begin{center}
\caption{Robust comparisons of several volatility estimation methods\label{tab45}}
\doublerulesep 0.5pt \small
\begin{tabular}{@{}|c| c|c|c|c|c|c|@{}} \hline \hline
 Measure & Empirical Formula & Hist & RiskM & Semi & NonBay &Integ\\
\hline\hline
IMADE     & Ave $(\times 10^{-6})$ & 0.5579 & 0.3025 & 0.4374 & 0.2748 &  0.2748\\
       & Relative Loss (\%)      & 103.01  & 10.08  & 59.17 &  0   & 0\\
\hline
MADE   & Ave $(\times 10^{-5})$ &  0.1115  & 0.1107 & 0.1111 & 0.1097 &   0.1061\\
       & Relative Loss (\%)     & 5.07  & 4.30 &  4.67 & 3.42 & 0\\
\hline
RADE   & Ave  ($\times 10^{-3}$)  &  0.4268  & 0.3901 & 0.4028 & 0.3885 & 0.3836\\
      & Relative Loss (\%)     & 11.27   & 1.71  & 5.00  & 1.28   & 0\\
\hline
ER     & Ave                    &  0.0628  & 0.0521 &  0.0493 &  0.0494  &  0.0428\\
\hline \hline
\end{tabular}
\end{center}
\end{table}

Table \ref{tab4} summarizes the results.
The historical simulation
approach has the smallest MADE, but suffers from poor forecast in
terms of IMADE. This is surprising. Why is it so different between IMADE and MADE?
This phenomenon may be produced by the non-stationarity of
the process. For the integrated method, even though the true volatility structure
is well captured because of the lowest
IMADE, extreme values of observations make the MADE quite large.
To more accurately calibrate the performance of the volatility estimation,
we use the $95\%$ up-trimmed mean instead of the mean to
summarize the values of the measures.
Table \ref{tab45} reports the trimmed means and the relative losses for different measures.
The similar conclusions
to those in Example 1 can be drawn from the table. This shows that our
integrated
method continues to perform better than other for this non-stationary case.
The Bayesian estimator performs comparably with the dynamically
integrated method and outperforms all others.

{\centering \subsection{Empirical Study}}

In this
section, we will apply the integrated volatility estimation
methods and others to the analysis of real financial data.

{\centering\subsubsection{Treasury Bond} }

We consider here the weekly returns of three treasury bonds
with terms 1, 5 and 10 years, respectively.

 We set the observations from January 4, 1974 to
December 30, 1994 as in-sample data, and those from January
6, 1995 up to August 8, 2003 as out-sample data.
The total sample size is $1545$ and the
in-sample size is $1096$.
The results are reported in Table \ref{tab5}.


\begin{table}[htbp]
\begin{center}  
\caption{Comparisons of several volatility estimation methods\label{tab5}}
\doublerulesep 0.5pt
\begin{tabular}{@{}|l|l|c|c|c|c|c|@{}} \hline \hline
Term  &  Measure    & Hist &  RiskM & Semi & NonBay & Integ \\
\hline
      & MADE & 0.01044 & 0.00787 & 0.00787 & 0.00794 &  0.00732 \\
1 year  & RADE & 0.05257 & 0.04231 & 0.04256 & 0.04225 & 0.04107\\
      & ER   & 0.022   & 0.020   & 0.022   & 0.016   & 0.038\\
\hline
      & MADE   & 0.01207 & 0.01253 & 0.01296 & 0.01278 & 0.01201\\
5 years      & RADE   & 0.05315  & 0.05494 & 0.05630 & 0.05562 & 0.05572\\
      & ER    & 0.007 & 0.014    & 0.016 &  0.011 & 0.058\\
\hline
      & MADE &  0.01041 & 0.01093 & 0.01103 & 0.01112 & 0.01018\\
 10 years     & RADE & 0.04939  & 0.05235 & 0.05296 & 0.05280 & 0.05151\\
      & ER   &  0.011  & 0.016 & 0.018   & 0.013 & 0.049\\
  \hline \hline
\end{tabular}
\end{center}
\end{table}

From Table \ref{tab5}, the integrated estimator is of the smallest
MADE and almost the smallest RADE, which reflects that the
integrated estimation method of the volatility is the best among
the five methods. Relative losses in MADE of the other estimators
with respect to the integrated estimator can easily be computed as
ranging from $8.47\%$
(NonBay) to $42.6\%$ (Hist) for the bond with one year term. For
the bonds with 5 or 10 years term, the five estimators
have close MADEs and RADEs, where the historical simulation method
is better than the RiskMetrics in terms of MADE and RADE, and the
integrated estimation approach has the smallest MADEs.
This demonstrates the advantage of using state
domain information which can help the time-domain prediction of
the changes in bond interest dynamics.

{\centering\subsubsection{Exchange Rate} }

We analyse the daily exchange rate of several foreign currencies
with US dollar. The data are from January 3, 1994 to August 1,
2003. The in-sample data consists of the observations before
January 1, 2001, and the out-sample data consists of the remaining
observations. The results are reported in Table \ref{tab6}. It is
seen that the integrated estimator has the smallest MADEs for the
exchange rates, which again supports our integrated volatility
estimation.\\

\begin{table}[htbp]
\begin{center}  
\caption{Comparisons of several volatility estimation methods\label{tab6}}
\doublerulesep 0.5pt
\begin{tabular}{@{}|l|l|c|c|c|c|c|@{}} \hline \hline
Currency  & Measure & Hist &  RiskM  & Semi & NonBay & Integ \\
\hline
        & MADE($\times 10^{-4}$) & 0.614  & 0.519  & 0.536  & 0.519  & 0.492\\
   U.K.   & RADE($\times 10^{-3}$) & 3.991 & 3.424 & 3.513 & 3.438 &3.491\\
        & ER   & 0.011    &   0.017   & 0.019 & 0.015 & 0.039\\
\hline
        & MADE($\times 10^{-4}$) & 0.172  & 0.132  & 0.135  & 0.135  & 0.126\\
 Australia       & RADE($\times 10^{-3}$) & 1.986 & 1.775 & 1.830 & 1.797 & 1.762\\
        & ER   & 0.054    & 0.025 & 0.026 & 0.022 & 0.043\\
\hline
        & MADE($\times 10^{-1}$) & 5.554    & 5.232   & 5.444 &  5.439   & 5.067\\
   Japan   & RADE($\times 10^{-1}$) &  3.596 & 3.546 & 3.622 & 3.588 & 3.560\\
        &  ER  & 0.014  & 0.011 & 0.019 & 0.012 & 0.029\\
\hline \hline
\end{tabular}
\end{center}
\end{table}

{\centering
\section{Conclusions}\label{sec6}}

We have proposed a Bayesian method and a dynamically integrated
method to aggregate the information from  the time-domain and the
state domain. The performance comparisons are studied both
empirically and theoretically. We have shown that the proposed
integrated method is effectively aggregating the information from
both the time and the state domains, and has advantages over some
previous methods.  It is powerful in forecasting volatilities for
the yields of bonds and for exchange rates.  Our study has also
revealed that proper use of information from both the time domain
and the state domain makes volatility forecasting more accurately.
Our method exploits the continuity in the time-domain and
stationarity in the state-domain.  It can be applies to situations
where these two conditions hold approximately.


{\centering\section{Appendix}\label{ap} }

\setcounter{equation}{0}
\renewcommand{\theequation}{A\arabic{equation}}

We collect technical conditions for the proof of our results.

\begin{description}
\item[(A1)] $\sigma^2(x)$ is Lipschitz continuous.

\item[(A2)] There exists a constant $L>0$ such that
$E|\mu(r_s)|^{2(p+\delta)}\leq L$ and
$E|\sigma(r_s)|^{2(p+\delta)}\leq L$ for any $s\in [t-\eta,t]$,
where $\eta$ is some positive constant, $p$ is an integer not less
than $4$ and $\delta>0$.

\item[(A3)] The discrete observations $\{r_{t_i} \}_{i=0}^N$ satisfy the
stationarity conditions of Banon (1978).  Furthermore, the $G_2$
condition of Rosenblatt (1970) holds for the transition operator.

\item[(A4)] The conditional density $p_\ell(y|x)$ of
$r_{t_{i+\ell}}$ given $r_{t_i}$ is continuous in the arguments $(y,x)$
and is bounded by a constant independent of $\ell$.

\item[(A5)] The kernel $W$ is a bounded, symmetric
probability density function with compact support, $[-1,1]$ say.

\item [(A6)] $(N - n) h \rightarrow \infty$, $(N- n) h^5 \to 0$, $(N-n) h \Delta
\to 0$.
\end{description}

Throughout the proof, we denote by $M$ a generic positive
constant, and use $\mu_s$ and $\sigma_s$ to represent $\mu(r_s)$
and $\sigma(r_s)$, respectively.

\noindent {\bf Proof  of Proposition \ref{P1}}. It suffices to
show that the process $\{r_s\}$ is H\"{o}lder-continuous with
order $q=(p-1)/(2p)$ and coefficient $K_1$, where
$E[K_1^{2(p+\delta)}]<\infty$, because this together with
assumption $(A1)$ gives the result of the lemma. By Jensen's
inequality and martingale moment inequalities (Karatzas \& Shreve
1991, Section 3.3.D), we have
\begin{align*}
E|r_u-r_s|^{2(p+\delta)}&\leq
M\left(E\left|\int_s^u\mu_vdv\right|^{2(p+\delta)}+
E\left|\int_s^u\sigma_vdW_v\right|^{2(p+\delta)}
\right)\\
&\leq M(u-s)^{2(p+\delta)-1}\int_s^uE|\mu_v)|^{2(p+\delta)}dv+
M(u-s)^{p+\delta-1}\int_s^uE|\sigma_v|^{2(p+\delta)}dv  \\
&\leq M(u-s)^{p+\delta}.
\end{align*}
Then by the Kolmogorov continuity
theorem (Revuz \& Yor 1991, Theorem 2.1),
$\{r_s\}$ is H\"{o}lder-continuous.
\\

\noindent{\bf Proof of Theorem \ref{T1}}. Let
$Z_{i,s}=(r_s-r_{t_i})^2$. Applying It\^{o} formula to $Z_{i,s}$,
we obtain
\begin{align*}
dZ_{i,s}=&2\Bigl(\int_{t_i}^s\mu_udu+\int_{t_i}^s\sigma_udW_u\Bigr)\Bigl(\mu_sds+\sigma_sdW_s\Bigr)+\sigma_s^2ds\\
=&2\left[\Bigl(\int_{t_i}^s\mu_udu+\int_{t_i}^s\sigma_udW_u\Bigr)\mu_sds+\sigma_s\Bigl(\int_{t_i}^s\mu_udu\Bigr)dW_s
\right]\\
&+2\Bigl(\int_{t_i}^s\sigma_udW_u\Bigr)\sigma_sdW_s+\sigma_s^2ds.
\end{align*}
Then $Y_i^2 $ can be decomposed as
$$Y_i^2=2a_i+2b_i+\bar{\sigma}_i^2, $$
where
$$a_i=\Delta^{-1}\left[\int_{t_i}^{t_{i+1}}\mu_sds\int_{t_i}^s\mu_udu+
\int_{t_i}^{t_{i+1}}\mu_sds\int_{t_i}^s\sigma_udW_u+\int_{t_i}^{t_{i+1}}
\sigma_sdW_s\int_{t_i}^s\mu_udu\right],$$
$$b_i=\Delta^{-1}\int_{t_i}^{t_{i+1}}\int_{t_i}^s\sigma_udW_u\sigma_sdW_s,$$
and
$$\bar{\sigma}_i^2=\Delta^{-1}\int_{t_i}^{t_{i+1}}\sigma_s^2ds.$$
Therefore, $\hat{\sigma}_{ES, t}^2
$ can be written as
\begin{eqnarray*}\label{019}
\hat{\sigma}_{ES,t}^2
&=&2\frac{1-\lambda}{1-\lambda^n}\sum_{i=t-n}^{t-1}\lambda^{t-i-1}a_i+
2\frac{1-\lambda}{1-\lambda^n}\sum_{i=t-n}^{t-1}\lambda^{t-i-1}b_i+
\frac{1-\lambda}{1-\lambda^n}\sum_{i=t-n}^{t-1}\lambda^{t-i-1}\bar{\sigma}_i^2\\
&\equiv& A_{n,\Delta} +B_{n,\Delta} +C_{n,\Delta}.
\end{eqnarray*}
By Proposition 1, as $n\Delta \rightarrow 0$,
$$
  |C_{n,\Delta}-\sigma^2_t |\leq K(n\Delta)^{q},
$$
where $q = (p-1)/(2p)$. This combined with Lemmas
\ref{L1}-\ref{L2} below completes the proof of the theorem.

\begin{lemma}\label{L1} If condition (A2) is satisfied, then
$E[A_{n,\Delta}^2]=O(\Delta).$
\end{lemma}

\noindent{\bf Proof }. Simple algrbea gives the result. In fact,
\begin{eqnarray*}
E(a_i^2) &\leq&
3E\left[\Delta^{-1}\int_{t_i}^{t_{i+1}}\mu_sds\int_{t_i}^s\mu_udu\right]^2
+3E\left[\Delta^{-1}\int_{t_i}^{t_{i+1}}\mu_sds\int_{t_i}^s\sigma_udW_u \right]^2 \\
&&+3E\left[\Delta^{-1}\int_{t_i}^{t_{i+1}}\sigma_sdW_s\int_{t_i}^s\mu_udu\right]^2\\
&\equiv& I_1(\Delta)+I_2(\Delta)+I_3(\Delta).
\end{eqnarray*}
Applying Jensen's inequality, we obtain that
\begin{eqnarray*}
I_1(\Delta)
&=&O(\Delta^{-1})E\Bigl[\int_{t_i}^{t_{i+1}}\int_{t_i}^s\mu_s^2\mu_u^2\,du\,ds\Bigr]\\
  &=& O(\Delta^{-1}) \int_{t_i}^{t_{i+1}}\int_{t_i}^sE(\mu_u^4+\mu_s^4)\,du\,ds
= O(\Delta).
\end{eqnarray*}
By Jensen's inequality, H\"{o}lder's inequality and martingale
moments inequalities, we have
\begin{eqnarray*}
I_2(\Delta)
 &=&O(\Delta^{-1})\int_{t_i}^{t_{i+1}}E\Bigl(\mu_s\int_{t_i}^s\sigma_u^2dW_u\Bigr)^2ds\\
  &=&O(\Delta^{-1})\int_{t_i}^{t_{i+1}}\Bigl\{E\Bigl[\mu_s\Bigr]^4
  E\Bigl[\int_{t_i}^{t_{i+1}}\sigma_udW_u \Bigr]^4\Bigr\}^{1/2}ds
  = O(\Delta).
\end{eqnarray*}
Similarly, $I_3(\Delta)=O(\Delta).$ Therefore,
$E(a_i^2)=O(\Delta)$. Then by the Cauchy-Schwartz inequality and
noting that $n(1-\lambda)=O(1)$, we obtain that
$$E[A^2_{n,\Delta}]
\leq
n\Bigl(\frac{1-\lambda}{1-\lambda^n}\Bigr)^2\sum_{i=1}^n\lambda^{2(n-i)}E(a_i^2)=O(\Delta)
.
$$

\begin{lemma}\label{L2} Under condition (A2), if $n\rightarrow \infty$ and $n\Delta\rightarrow
0$, then
\begin{equation}\label{028}
s_{1,t}^{-1}\sqrt{n}B_{n,\Delta}\stackrel{\mathcal
D}{\longrightarrow} {\mathcal N}\Bigl(0,1\Bigr).
\end{equation}
\end{lemma}

\noindent{\bf Proof}. Note that
\[b_j=\sigma_{t}^2\Delta^{-1}\int_{t_j}^{t_{j+1}}(W_s-W_{t_j})dW_s+\epsilon_j, \]
where
\[\epsilon_j=\Delta^{-1}\int_{t_j}^{t_{j+1}}(\sigma_s-\sigma_{t})
\left[\int_{t_j}^s\sigma_udW_u\right]dW_s+
\Delta^{-1}\sigma_{t}\int_{t_j}^{t_{j+1}}
\left[\int_{t_j}^s(\sigma_u-\sigma_{t})dW_u\right]dW_s.
\]
By the central limit theorem for martingale (see Hall \& Heyde
1980, Corollary 3.1), it suffices to show that
\begin{equation}\label{007}
V_n^2\equiv E[s_{1,t}^{-2}{n}B_{n,\Delta}^2] \rightarrow 1,
\end{equation}
and the following Lyapunov condition holds:
\begin{equation}\label{008}
\sum_{i=t-n}^{t-1}E\left(\sqrt{n}\frac{1-\lambda}{1-\lambda^n}\lambda^{t-i-1}b_i\right)^4\rightarrow0.
\end{equation}
Note that
\begin{eqnarray}
\frac{\Delta^2}{2}E(\epsilon_j^2) &\leq&
E\Bigl\{\int_{t_j}^{t_{j+1}}
(\sigma_s-\sigma_{t})\Bigl[\int_{t_j}^s\sigma_udW_u\Bigr]dW_s\Bigr\}^2\nonumber\\
&&+ \sigma_t^2E\Bigl\{\int_{t_j}^{t_{j+1}}
\Bigl[\int_{t_j}^s(\sigma_u-\sigma_{t_t})dW_u\Bigr]dW_s\Bigr\}^2\nonumber\\
&\equiv& L_{n1}+L_{n2}.\label{024}
\end{eqnarray}
 By Jensen's inequality, H\"{o}lder's inequality and moments
inequalities for martingale, we have
\begin{eqnarray}
L_{n1} &\leq& \int_{t_j}^{t_{j+1}}
E\Bigl\{(\sigma_s-\sigma_{t})^2\Bigl[\int_{t_j}^s\sigma_udW_u\Bigr]^2\Bigr\}\,ds\nonumber\\
&\leq& \int_{t_j}^{t_{j+1}}
\Bigl\{E(\sigma_s-\sigma_{t})^4E\Bigl[\int_{t_j}^s\sigma_udW_u\Bigr]^4\Bigr\}^{1/2}\,ds\nonumber\\
&\leq& \int_{t_j}^{t_{j+1}}
\Bigl\{E[K(n\Delta)^q]^4 \,36\Delta\int_{t_j}^s E(\sigma_u^4)du\Bigr\}^{1/2}\,ds\nonumber\\
&\leq& M(n\Delta)^{2q}\Delta^2. \label{010}
\end{eqnarray}
Similarly,
\begin{equation}\label{011}
L_{n2}\leq
 M(n\Delta)^{2q}\Delta^2.
\end{equation}
By (\ref{024}), (\ref{010}) and (\ref{011}),
\begin{equation}\label{027}
E(\epsilon_j^2)\leq M(n\Delta)^{2q} .
\end{equation}
Therefore,
$$E[\sigma^{-4}_tb_j^2]=\frac{1}{2}+O((n\Delta)^q).$$
By the theory of stochastic calculus, simple algebra gives that
$E(b_j)=0$ and $E(b_ib_j)=0$ for $ i\neq j$. It follows that
$$V_n^2=E(s_{1,t}^{-2}{n}B_{n,\Delta}^2)=\sum_{i=t-n}^{t-1}
E\left(2s_{1,t}\sqrt{n}\frac{1-\lambda}{1-\lambda^n}\lambda^{t-i-1}b_i\right)^2
\rightarrow 1.
$$
That is, (\ref{007}) holds. For (\ref{008}), it suffices to prove
that $E(b_j^4)$ is bounded, which holds by applying the moment
inequalities for martingales to $b_j^4$.

\noindent{\bf Proof of Theorem \ref{T2}}. The proof is completed
by using the same lines in Fan \& Zhang (2003).

\noindent{\bf Proof of Theorem \ref{T3}}. By Fan \& Yao (1998),
the volatility estimator $\hat{\sigma}^2_{S,t_N}$ behaves as if
the instantaneous return function $f$ is known, hence without loss
of generality we assume that $f(x)=0$ and hence $\hat{R}_i=Y_i^2.$
Let $\mathbf{Y}=(Y_{0}^2,\cdots,Y_{N-n-1}^2)^T$,
$\mathbf{W}=\text{diag}\{W_h(r_{t_0}-r_{t_N}),\cdots,W_h(r_{t_{N-n-1}}-r_{t_N})\},$
and
$$\mathbf{X}=\left(\begin{array}{cc}
1 & r_{t_0}-r_{t_N}\\
\vdots & \vdots\\
1 & r_{t_{N-n-1}}-r_{t_N}
\end{array}\right). $$
Denote by $m_{i}=E[Y_i^2|r_{t_i}]$,
$\mathbf{m}=(m_{0},\cdots,m_{N-n-1})^T$ and $\mathbf{e_1}=(1,0)^T$.
Define $\mathbf{S_N=X^TWX}$ and $\mathbf{T_N=X^TWY}$. Then it can
be written that (see Fan \& Yao, 2003)
$$\hat{\sigma}_{S,t_N}^2=\mathbf{e_1^TS_N^{-1}T_N}.$$
Hence
\begin{eqnarray}
\hat{\sigma}_{S,t_N}^2-\sigma^2_{t_N}
&=&\mathbf{e_1^TS_N^{-1}X^TW} \{\mathbf{m}-\mathbf{X}\bbeta_N\}
+\mathbf{e_1^TS_N^{-1}X^TW(Y-m)}\nonumber\\
&\equiv&\mathbf{e_1^Tb}+\mathbf{e_1^Tt},
\end{eqnarray}
where $\bbeta_N=(m(r_{t_N}),\ m'(r_{t_N}))^T$ with
$m(r_{t_N})=E[Y_1^2|r_{t_1}=r_{t_N}]$. By Fan \& Zhang (2003), the
bias vector $\mathbf{b}$ converges in probability to a vector
$\mathbf{\bar{b}}$ with
$\mathbf{\bar{b}}=O(h^2)=o(1/\sqrt{(N-n)h})$. In the following, we
will show that the centralized vector $\mathbf{t}$ is
asymptotically normal.

In fact, put
$\mathbf{u}=(N-n)^{-1}\mathbf{H^{-1}X^TW(Y-m)}$ where
$\mathbf{H}=\text{diag}\{1,h\}$,
then by Fan \& Zhang (2003)
the vector $\mathbf{t}$ can
be written as
\begin{equation}\label{032}
\mathbf{t}=p^{-1}(r_{t_N})\mathbf{H^{-1}S^{-1}u}(1+o_p(1)),
\end{equation}
where $\mathbf{S}=(\mu_{i+j-2})_{i,j=1,2}$ with $\mu_j=\int
u^jW(u)du$. For any constant vector $\mathbf{c}$, define
$$Q_N=\mathbf{c^Tu}=\frac{1}{N-n}\sum_{i=0}^{N-n-1}\{Y_i^2-m_i\}C_h(r_{t_i}-r_{t_N}),$$
where $C_h(\cdot)=1/h C(\cdot/h)$ with $C(x)=c_0W(x)+c_1xW(x)$. Applying the ``big-block" and
``small-block" arguments in Fan \& Yao (2003, Theorem 6.3), we
obtain
\begin{equation}\label{029}
 \theta^{-1}(r_{t_N}) \sqrt{(N-n)h}Q_N \toD N\left(0,1\right),
\end{equation}
where
$\theta^2(r_{t_N})=2p(r_{t_N})\sigma^4(r_{t_N})\int_{-\infty}^{+\infty}C^2(u)du$.
In the following, we will decompose $Q_N$ into two parts, $Q_N'$
and $Q_N''$, which satisfy that
\begin{itemize}
\item[(i)]  $(N-n)hE[\theta^{-1}(r_{t_N})Q_N']^2\leq\frac{h}{N-n}\left(h^{-1}a_N(1+o(1))+
(N-n)o(h^{-1}) \right)\rightarrow 0.$
\item[(ii)] $Q_N''$ is identically distributed as $Q_N$
and is asymptotically independent of $\hat{\sigma}^2_{ES,t_N}$.
\end{itemize}
 Define
\begin{equation}\label{009}
Q_N'=\frac{1}{N-n}\sum_{i=0}^{a_N}\{Y_i^2-E[Y^2_i|r_{t_i}]\}C_h(r_{t_i}-r_{t_N}),
\end{equation}
and $$Q_N''=Q_N-Q_N',$$ where $a_N$ is a positive integer
satisfying $a_N=o(N-n)$ and $a_N\Delta\rightarrow\infty$. Let
$\vartheta_{N,\ell}=(Y_i^2-m_i)C_h(r_{t_i}-r_{t_N})$, then by Fan
\& Zhang (2003)
\begin{equation}\label{002}
\Var[\theta^{-1}(r_{t_N})\vartheta_{N,1}]=h^{-1}(1+o(1)) \text{
and }
\sum_{\ell=1}^{N-n-2}|\Cov(\vartheta_{N,1},\vartheta_{N,\ell+1})|=
o(h^{-1}),
\end{equation}
which yields the result in (i).
This combined with (\ref{029}), (i) and (\ref{009}) leads
to
\begin{equation}\label{030}
\theta^{-1}(r_{t_N}) \sqrt{(N-n)h}Q_N'' \toD N\left(0,1\right).
\end{equation}
Note that the stationarity conditions of Banon (1978) and the $G_2$
condition of Rosenblatt (1970) on the transition operator imply
that the $\rho$-mixing coefficient $\rho(\ell)$ of $\{r_{t_i}\}$ decays exponentially,
 and the strong-mixing
coefficient $\alpha(\ell)\leq\rho(\ell)$, it follows that
\begin{equation}
\left|E\exp\{i\xi(Q_N''+\hat{\sigma}^2_{ES,t_N})\}-E\exp\{i\xi(Q_N''\}
E\exp\{i\xi\hat{\sigma}^2_{ES,t_N})\}\right|\leq
32\alpha(s_N)\rightarrow 0,
\end{equation}
for any $\xi\in\mathbb{R}$. Using the theorem of Volkonskii \&
Rozanov (1959), one gets the asymptotic independence of
$\hat{\sigma}^2_{ES,t_N}$ and $Q_N''$.

By (i), $\sqrt{(N-n)h}Q_N'$ is
asymptotically negligible.  This together with Theorem \ref{T1}
lead to
$$d_1\theta^{-1}(r_{t_N})\sqrt{(N-n)h}Q_N+d_2V_2^{-1/2}\sqrt{n}[\hat{\sigma}^2_{ES,t_N}-\sigma^2(r_{t_N})]
\stackrel{\mathcal D}{\longrightarrow}{\mathcal N}
\Bigl(0, d_1^2+d_2^2\Bigr),
$$
for any $d_1,\ d_2\in\mathbb{R}$, where
$V_2=\frac{e^c+1}{e^c-1}\sigma^4(r_{t_N})$. Since $Q_N$ is a
linear transform of $\mathbf{u}$,
\begin{eqnarray*}
{\mathbf V}^{-1/2}\left[ \begin{array}{c}  \sqrt{(N-n)h}\mathbf{u} \\
                          \sqrt{n}[\hat{\sigma}^2_{ES,t_N}-\sigma^2(r_{t_N})]
           \end{array} \right]
\stackrel{\mathcal D}{\longrightarrow} {\mathcal N}(0, I_{3}),
\end{eqnarray*}
where ${\mathbf V}=\mbox{\rm blockdiag}\{V_1,V_2\}$ with $V_1=
2\sigma^4(r_{t_N})p(r_{t_N})\mathbf{S^*}$,
where $\mathbf{S^*}=(\nu_{i+j-2})_{i,j=1,2}$ with $\nu_j=\int
u^jW^2(u)du$. This combined with (\ref{032}) gives the joint
asymptotic normality of $\mathbf{t}$ and
$\hat{\sigma}^2_{ES,t_N}$. Note that
$\mathbf{b}=o_p(1/\sqrt{(N-n)h})$, it follows that
\begin{eqnarray*}
\Sigma^{-1/2}\Bigl( \begin{array}{c}  \sqrt{(N-n)h}[\hat{\sigma}^2_{S,t_N}-\sigma^2(r_{t_N})] \\
                          \sqrt{n}[\hat{\sigma}^2_{ES,t_N}-\sigma^2(r_{t_N})]
           \end{array} \Bigr)
\stackrel{\mathcal D}{\longrightarrow} {\mathcal N}(0,I_{2}),
\end{eqnarray*}
where $\Sigma=\mbox{\rm
diag}\{2\sigma^4(r_{t_N})\nu_0/p(r_{t_N}),V_2\}$. Note that
$\hat{\sigma}^2_{S,t_N}$ and $\hat{\sigma}^2_{ES,t_N}$ are
asymptotically independent, it follows that the asymptotical
normality of $\hat{\sigma}_{I,t_N}^2$ holds.

\begin{center}
\textbf{Acknowledgements}
\end{center}

 The work was partially supported by a grant from
the Research Grants Council of the Hong Kong SAR (Project No. CUHK
400903/03P), the NSF grant DMS-0355179 and the Chinese NSF grants
10471006 and 10001004.
The authors thank Dr. Juan Gu for various assistances.

\hskip2mm
\bigskip
\begin{center}
{\bf References}
\end{center}
\begin{singlespace}
\begin{itemize}
\item[] Banon, G. (1978). Nonparametric identification for diffusion
processes. {\em SIAM J. Control Optim} {\bf 16}, 380-395.

\item[] Barndoff-Neilsen, O.E. \& Shephard, N. (2001).
    Non-Gaussian
    Ornstein-Uhlenbeck-based models and some of their uses
    in financial economics (with discussion.
    {\em J. R. Statist. Soc. B} {\bf 63}, 167-241.

\item[] Barndoff-Neilsen, O.E. \& Shephard, N. (2002).  Econometric
    analysis of realized volatility and its use in
    estimating stochastic volatility models.
    {\em J. R. Statist. Soc. B} {\bf 64}, 253-280.

\item[] Bollerslev, T. \& Zhou, H. (2002).  Estimating stochastic
        volatility diffusion using conditional moments of
        integrated volatility.  {\em Jour. Econometrics} {\bf 109},
        33-65.

\item[] Chan, K.C., Karolyi, A.G., Longstaff, F.A. \& Sanders, A.B.
        (1992).  An empirical comparison of alternative models of the
        short-term interest rate.  {\it Journal of Finance} {\bf 47},
        1209-1227.
\item[] Chapman, D.A. \& Pearson, N.D. (2000).
      Is the short rate drift actually nonlinear?
      {\it Journal of Finance} {\bf 55}, 355--388.

\item[] Cox, J.C., Ingersoll, J.E. \& Ross, S. A. (1985).  A theory of
    the term structure of interest rates.  {\it Econometrica} {\bf 53},
        385-467.

\item[] Dav\'e, R. D. \& Stahl, G. (1997).
             On the accuracy of VaR estimates based on
         the Variance-Covariance approach. Working paper,
         Olshen \& Associates.

\item[] Duffie, D. \& Pan, J. (1997). An overview of Value at Risk.
      {\it The Journal of Derivatives}, 7--49.



\item[] Fan, J. \& Gu, J. (2003). Semiparametric estimation of
value-at-risk.   {\it Econometrics Journal} {\bf 6}, 261-290.



\item[] Fan, J. \& Yao, Q. (1998).
      Efficient estimation of conditional variance functions
      in stochastic regression. {\it Biometrika} {\bf 85}, 645--660.

\item[] Fan, J. \& Yao, Q. (2003).
{\em Nonlinear Time Series: Nonparametric and Parametric Methods}, Springer-Verlag, New York.

\item[] Fan, J. \& Zhang, C.M. (2003). A Reexamination of
      Diffusion Estimators with Applications to Financial Model Validation.
      {\em J. Am. Statist. Assoc.} {\bf 98}, 118--134.



\item[] Genon-Catalot, Jeanthheau, T. \& Laredo, C. (1999). {
Parameter estimation for discretely observed stochastic volatility
models}, {\it Bernoulli} {\bf 5}, 855--872.

\item[] Gijbels, I., Pope, A., \& Wand, M.P. (1999). Understanding
    exponential smoothing via kernel regression. {\sl J.
    R. Statist. Soc. B} {\bf 61}, 39--50.


\item[] Hall, P. \& Heyde, C. (1980). {\em Martingale limit theorem
and its applications.}
       Academic Press.





\item[] Kloeden, D.E., Platen, E., Schurz, H. \& S\o rensen, M. (1996).
      On effects of discretization on estimators of drift parameters for
      diffusion processes. {\it Journal of Applied Probability}
     {\bf 33}, 1061--1076.





\item[] Karatzas, I. \& Shreve, S. (1991). {\em Brownian motion and
stochastic calculus (2nd edition).} Springer-Verlag, New York.

\item[] Morgan, J.P. (1996) {\it RiskMetrics Technical Document},
      Fourth edition, New York.



\item[] Stanton, R. (1997). A nonparametric models of term structure dynamics
     and the market price of interest rate risk.
      {\it Journal of Finance} {\bf LII}, 1973--2002.


\item[] Revuz, D. \& Yor, M. (1991). {\em Continuous Martingales
and Brownian Motion.}  Springer-Verlag.

\item[] Rosenblatt, M. (1970). Density estimates and Markov
sequences. In {\em Nonparametric Techniques in Statistical
Inference }(ML Puri, ed.) 199-213. Cambridge Univ. Press.

\item[] Ruppert, D., Wand, M.P., Holst, U. \& H\"ossjer, O. (1997).
      Local polynomial variance function estimation.  {\em Technometrics}
      {\bf 39}, 262-273.

\item[] Spokoiny, V. (2000). Drift estimation for nonparametric diffusion.
{\it Ann. Statist.} {\bf 28}, 815--836.

\item[] Tong, H. (1990).
      {\it Non-Linear Time Series: A Dynamical System Approach}.
      Oxford University Press, Oxford.

\item[] Tong, H. (1995). A personal overview of non-linear time series
        analysis from a chaos perspective (with discussion).
        {\sl Scandinavian Journal of Statistics} {\bf 22}, 399-445.

\item[] Zhang, C.M. (2003). Calibrating the degrees of freedom for automatic
data-smoothing and effective curve-checking.
{\it J. Am. Statist. Assoc.} {\bf 98}, 609-628 .

\end{itemize}
\end{singlespace}
\end{document}